# О СЛЕДАХ СРЕДНИХ СПЕКТРАЛЬНЫХ РАЗЛОЖЕНИЙ, СООТВЕТСТВУЮЩИХ САМОСОПРЯЖЕННЫМ ЭЛЛИПТИЧЕСКИМ ПСЕВДОДИФФЕРЕНЦИАЛЬНЫМ ОПЕРАТОРАМ РАСПРЕДЕЛЕНИЙ В КЛАССАХ СОБОЛЕВА–ЛИУВИЛЛЯ И НИКОЛЬСКОГО-БЕСОВА

Бабаев М.М[1]

Пусть $\Omega$ – произвольная $N$-мерная ограниченная область с гладкой границей класса $C^\infty$, либо $\Omega = R^N$.

Определим класс символов $S^m(\Omega)$ как множество функций $a(x,y) \in C^\infty\left(\Omega \times \{R^N \setminus 0\}\right)$ таких, что для любого компакта $K \Subset \Omega$ и любых мультииндексов $\alpha, \beta$ выполняется неравенство

$$\left|D_x^\alpha D_y^\beta a(x,y)\right| \leq const \cdot (1+|y|)^{m-|\beta|} \qquad (1)$$

с постоянной, не зависящей от $x \in K$ и $y \in R^N \setminus 0$.

Каждому символу $a(x,y) \in S^m(\Omega)$ формула

$$A(x,D)f(x) = \int_{R^N} a(x,y) \hat{f}(y) \exp(ixy) dy, \qquad (2)$$

сопоставляет псевдодифференциальный оператор $A(x,D) \in OPS^m(\Omega)$, где

$$\hat{f}(y) = \frac{1}{(2\pi)^N} \int_{R^N} f(x) \exp(-ixy) dx. \qquad (3)$$

Пусть $A$ – какое-нибудь положительное самосопряженное расширение в $L_2(\Omega)$ эллиптического оператора $A(x,D) \in OPS^m(\Omega)$ со скалярными однородными символами степени $m \geq 1$ и с постоянными коэффициентами, где $\Omega$ – произвольная $N$-мерная ограниченная область с гладкой границей

[1] **Бабаев Махкамбек Мадаминович** Национальний Университет Узбекистана, г. Ташкент, Узбекистан, *babayevm@mail.ru*

класса $C^\infty$, либо $\Omega = R^N$ (определения этих классов символов см. в [1, с.45, 70]).

Определим оператор

$$p(tA)f(x) = \int_{-\infty}^{\infty} \hat{p}(z)\exp(iztA)f(x)dz, \quad (4)$$

где $\hat{p}(z) = \dfrac{1}{2\pi} \int_{-\infty}^{\infty} p(\lambda)\exp(-i\lambda z)d\lambda$, причем оператор $\exp(iztA)$ понимается как в спектральной теореме Дж. фон Неймана, т.е.

$$\exp(iztA) = \int_{-\infty}^{\infty} \exp(izt\lambda)dE_\lambda, \quad (5)$$

где $E_\lambda$ – разложение единицы, порождаемое оператором $A$. Всюду через $\sigma_A = \sigma(y)$ будем обозначать символ оператора $A$. В частности, если положить

$$p(z) = \left(1-|z|\right)_+^s = \begin{cases} (1-|z|)^s & \text{при } |z| \leq 1, \\ 0 & \text{при } |z| > 1, \end{cases} \quad (6)$$

$t = \dfrac{1}{\lambda}$, то получим среднее Рисса $E_\lambda^s f(x) = p(\dfrac{1}{\lambda}A)f(x)$ спектрального разложения $E_\lambda f(x)$ порядка $s \geq 0$.

**Определение 1.** Пусть $s \in R$, $1 \leq p \leq \infty$. Положим

$$\|f\|_p^s = \left\| F^{-1}\left\{(1+|\xi|^2)^{\frac{s}{2}} Ff(\xi)\right\} \right\|_{L_p(R^N)}, \quad (7)$$

где $Ff(\xi) = \int_{R^N} \exp(-ix\xi)f(x)dx$ – преобразование Фурье. Пространство Лиувилля $L_p^s(R^N)$ определяются соотношениями

$$L_p^s(R^N) = \left\{f : f \in S'(R^N), \|f\|_p^s < \infty\right\}. \quad (8)$$

**Определение 2.** Пусть $\Omega \subset R^N$ — произвольная (ограниченная или неограниченная) область, $-\infty < s < +\infty$, $1 \leq p \leq \infty$. Пространства Лиувилля $L_p^s(\Omega)$ — это сужение пространства $L_p^s(R^N)$ на $\Omega$ и

$$\|f\|_{L_p^s(\Omega)} = \inf_{\substack{g|_\Omega = f \\ g \in L_p^s(R^N)}} \|g\|_{L_p^s(R^N)} . \qquad (9)$$

Обозначим символом $\overset{0}{L}{}_p^\alpha(\Omega)$ класс функций, принадлежащих пространству Лиувилля $L_p^\alpha(\Omega)$ и имеющих компактный в $\Omega$ носитель (определение классов $L_p^\alpha(R^N)$ см. в [2], [3], для многообразий эти классы, как обычно, определяются путем перехода к локальным координатам [4]).

Следуя Я. Петре [5]-[8], дадим определение пространства Никольского-Бесова. Пусть функция $\varphi(\xi) \in S(R^N)$ такая, что

$$\operatorname{supp} \varphi = \{\xi : 2^{-1} \leq |\xi| \leq 2\}, \quad \varphi(\xi) > 0 \text{ при } 2^{-1} \leq |\xi| \leq 2,$$

$$\sum_{k=-\infty}^{+\infty} \varphi\left(2^{-k}\xi\right) = 1 \text{ при } \xi \neq 0 . \qquad (10)$$

Определим функции $\varphi_k(\xi)$ и $\psi(\xi)$ соотношениями

$$F\varphi_k(\xi) = \varphi(2^{-k}\xi), \quad k = 0, \pm 1, \pm 2, \ldots ,$$

$$F\psi(\xi) = 1 - \sum_{k=1}^{\infty} \varphi(2^{-k}\xi) , \qquad (11)$$

где $Ff(\xi) = \int_{R^N} \exp(-ix\xi) f(x) dx$ - преобразование Фурье.

**Определение 3.** Пусть $s \in R$, $1 \leq p \leq \infty$ и $1 \leq q \leq \infty$. Положим

$$\|f\|_{p\,q}^{s} = \|\psi * f\|_p + \left(\sum_{k=1}^{+\infty} (2^{sk} \|\varphi_k * f\|_p)^q\right)^{\frac{1}{q}} \qquad (12) .$$

Пространства Никольского-Бесова $B_{pq}^s(R^N)$ определяются соотношениями

$$B_{pq}^s(R^N) = \left\{ f : f \in S'(R^N), \|f\|_{pq}^s < \infty \right\}. \qquad (13)$$

Следующая теорема доставляет другое определение пространств Никольского-Бесова $B_{pq}^s(R^N)$ $(s>0)$ в терминах производных и модулей непрерывности. Модуль непрерывности определяется формулой

$$\omega_p^m(t, f) = \sup_{|y|<t} \|\Delta_y^m f\|_p, \qquad (14)$$

где $\Delta_y^m$ – оператор взятия разности порядка $m$:

$$\Delta_y^m f(x) = \sum_{k=0}^{m} \binom{m}{k} (-1)^k f(x+ky) \qquad (15).$$

Справедлива следующая

**Теорема (Й. Берг, Й. Лёфстрём [2]).** Пусть $s>0$, а $m$ и $N_1$ – целые числа, такие, что $m+N_1>s$ и $0 \leq N_1 < s$. Тогда при $1 \leq p \leq \infty$, $1 \leq q \leq \infty$ справедлива эквивалентность нормы

$$\|f\|_{pq}^s \Box \|f\|_p + \sum_{j=1}^{N} \left( \int_0^\infty \left( t^{N_1-s} \, \omega_p^m(t, \frac{\partial^{N_1} f}{\partial x_j^{N_1}}) \right)^q \frac{dt}{t} \right)^{\frac{1}{q}}. \qquad (16)$$

Если $s>0$, $1 \leq p < \infty$ и $1 \leq q < \infty$, то пространства Бесова [3], [9] имеет вид

$$B_{pq}^s(R^N) = \left\{ f: \quad f \in W_p^{[s]^\lnot}(R^N), \quad \|f\|_{pq}^s = \|f\|_{W_p^{[s]^\lnot}(R^N)} + \right.$$

$$+ \sum_{|\alpha|=[s]^-} \left( \int_{R^N} |h|^{-\{s\}^+ q} \left\| \Delta_h^2 D^\alpha f \right\|_{L_p(R^N)}^q \frac{dh}{|h|^N} \right)^{\frac{1}{q}} < \infty \right\}, \qquad (17)$$

где $s = [s]^- + \{s\}^+$, $s = [s]^-$ — целое число и $0 < \{s\}^+ \leq 1$. Здесь если $1 \leq p < \infty$ и $m = 1, 2, 3, \ldots$, то

$$W_p^m(R^N) = \left\{ f : f \in L_p(R^N), \ \|f\|_{W_p^m(R^N)} = \sum_{|\alpha| \leq m} \left\| D^\alpha f \right\|_{L_p(R^N)} < \infty \right\} \qquad (18)$$

пространства Соболева [4], [10], [11], [12], [13]–[15]. А также определим при $1 \leq p < \infty$ и $s > 0$, $s$ — нецелое, то пространства Слободецкого [15]

$$W_p^s(R^N) = \left\{ f : f \in W_p^{[s]}(R^N), \ \|f\|_{W_p^s(R^N)} = \|f\|_{W_p^{[s]}(R^N)} + \right.$$

$$\left. + \sum_{|\alpha|=[s]} \left( \iint_{R^N \times R^N} \frac{|D^\alpha f(x) - D^\alpha f(y)|}{|x-y|^{N+\{s\}p}} dx dy \right)^{\frac{1}{p}} < \infty \right\}. \quad (19)$$

Если $s > 0$, $1 \leq p < \infty$, $q = \infty$, то пространства Никольского [4], [9] имеет вид

$$H_p^s(R^N) = B_{p\infty}^s(R^N) = \left\{ f : f \in W_p^{[s]^-}(R^N), \ \|f\|_{p\infty}^s = \|f\|_{W_p^{[s]^-}(R^N)} + \right.$$

$$\left. + \sum_{|\alpha|=[s]^-} \sup_{0 \neq h \in R^N} |h|^{-\{s\}^+} \left\| \Delta_h^2 D^\alpha f \right\|_{L_p(R^N)} < \infty \right\}. \quad (20)$$

**Определение 4.** Пусть $\Omega \subset R^N$ – произвольная (ограниченная или неограниченная) область, $-\infty < s < +\infty$, $1 \le p \le \infty$ и $1 \le q \le \infty$. Пространства $B^s_{p\,q}(\Omega)$ – это сужение пространства $B^s_{p\,q}(R^N)$ на $\Omega$,

$$\|f\|_{B^s_{p\,q}(\Omega)} = \inf_{\substack{g|_\Omega = f \\ g \in B^s_{p\,q}(R^N)}} \|g\|_{B^s_{p\,q}(R^N)}. \qquad (21)$$

Обозначим через $\overset{0}{B}{}^{\alpha}_{p\,q}(\Omega)$ класс функций, принадлежащих пространству $B^{\alpha}_{p\,q}(\Omega)$ и имеющих компактный в $\Omega$ носитель.

Одно из эквивалентных определений пространства Бесова $B^s_{p,q}(\Omega)$ дается с помощью вещественного метода интерполяции:

$$B^s_{p,q}(\Omega) = \left(L^{s_0}_p(\Omega),\ L^{s_1}_p(\Omega)\right)_{\theta,\,q}, \qquad (22)$$

где $1 \le p,\ q \le \infty$, $0 < \theta < 1$, $s = (1-\theta)s_0 + \theta s_1$ (см., например, [2], [14-15]). При $1 < p < \infty$ мы имеем также

$$W^s_p(R^N) = \begin{cases} L^s_p(R^N) & \text{при } s = 0, 1, 2, \ldots, \\ B^s_{p\,p}(R^N) & \text{при } 0 < s \ne 1, 2, \ldots \end{cases} \qquad (23)$$

Справедлива следующая

**Теорема 1.** Пусть функция $p(\lambda)$, определенная в $\overline{R}_+ = [0, \infty)$, удовлетворяет условиям

1) $\int\limits_0^\infty |p(\lambda)| \lambda^{\frac{N-\alpha_0}{m}-1} d\lambda < \infty$ и $p(0) = 1$;

2) $p(\lambda) \in C^l(\overline{R}_+)$ и $\left|p^{(j)}(\lambda)\right| \le C_j(1+\lambda)^{-j}$, $j = \overline{0, l}$, где

$l$ – неотрицательное целое число;

3) $l > N\left(\dfrac{1}{2} - \dfrac{1}{p_0}\right)$, $\alpha_0 = \dfrac{N}{p_0}$, $\alpha \geq 0$, $2 \leq p \leq p_0 < \infty$, $\varepsilon = N\left(\dfrac{1}{p} - \dfrac{1}{p_0}\right)$, $\beta \geq \alpha_0 + \alpha + \varepsilon$.

Тогда для любой непрерывной функции $u(x) \in \overset{0}{L}{}_p^\beta(\Omega)$ с компактным носителем выполняется равенство

$$\lim_{t \to +0} \left\| p(tA)u(x) - u(x) \right\|_{L_p^\alpha(M)} = 0 \tag{24}$$

для любого компакта $M \Subset \Omega$.

**Теорема 2.** Пусть функция $p(\lambda)$, определенная в $\overline{R}_+ = [0, \infty)$, удовлетворяет условиям $p(0) = 1$ и $p(\lambda) \in L_\infty(\overline{R}_+) \cap C([0, \tau))$, $\tau > 0$. Пусть $\alpha_0 > \dfrac{N}{p_0}$, $\alpha \geq 0$, $1 < p \leq p_0 \leq 2$ (либо $p = p_0 = 1$), $\varepsilon = N\left(\dfrac{1}{p} - \dfrac{1}{p_0}\right)$, $\beta \geq \alpha_0 + \alpha + \varepsilon$. Тогда для любой непрерывной функции $u(x) \in \overset{0}{L}{}_p^\beta(\Omega)$ с компактным носителем выполняется равенство

$$\lim_{t \to +0} \left\| p(tA)u(x) - u(x) \right\|_{L_p^\alpha(M)} = 0$$

для любого компакта $M \Subset \Omega$.

Доказательство теоремы 1 и 2 приведено в работе [16].

**Определение 5.** Средних спектральным разложением для распределения $f \in E'(\Omega)$ назовем распределения $p(tA)f \in D'(\Omega)$, действующее по правилу

$$< p(tA)f, \varphi > \; = \; < f, p(tA)\varphi > \tag{25}$$

для любой функции $\varphi \in D(\Omega) = C_0^\infty(\Omega)$.

Корректность данного определения следует из того, что для любого $t > 0$ оператор $p(tA)$ непрерывен из $D(\Omega)$ в $E(\Omega)$. Более того, справедлива следующее утверждение.

**Предложение 1.** Если выполнено условия теоремы 1 (или 2), то для любой функции $\varphi \in D(\Omega) = C_0^\infty(\Omega)$ в топологии $E(\Omega) = C^\infty(\Omega)$ выполняется равенство

$$\lim_{t \to +0} p(tA)\varphi(x) = \varphi(x).$$

Отсюда следует, что для любого распределения $f \in E'(\Omega)$ и любой функции $\varphi \in D(\Omega) = C_0^\infty(\Omega)$ выполняется равенство

$$\lim_{t \to +0} <p(tA)f, \varphi> \; = \; <f, \varphi>.$$

Иначе говоря, справедливо следующее утверждение.

**Предложение 2.** Если выполнено условия теоремы 1 (или 2), то для каждого распределения $f \in E'(\Omega)$ в топологии $D'(\Omega)$ имеет место соотношение

$$\lim_{t \to +0} p(tA)f = f. \qquad (26)$$

Более того, из теоремы 1 следует справедливость следующая

**Теорема 3.** Пусть функция $p(\lambda)$, определенная в $\overline{R}_+ = [0, \infty)$, удовлетворяет условиям

1) $\int_0^\infty |p(\lambda)| \lambda^{\frac{N-\alpha_0}{m}-1} d\lambda < \infty$ и $p(0) = 1$;

2) $p(\lambda) \in C^l(\overline{R}_+)$ и $|p^{(j)}(\lambda)| \leq C_j(1+\lambda)^{-j}$, $j = \overline{0, l}$, где

$l$ – неотрицательное целое число;

3) $l > N\left(\dfrac{1}{2} - \dfrac{1}{p_0}\right)$, $\alpha_0 = \dfrac{N}{p_0}$, $\alpha \geq 0$, $2 \leq p \leq p_0 < \infty$, $\varepsilon = N\left(\dfrac{1}{p} - \dfrac{1}{p_0}\right)$, $\beta \leq \alpha - \alpha_0 - \varepsilon$.

Тогда для любого распределения $f \in E'(\Omega)$ из класса $\overset{0}{L}{}_q^{-\beta}(\Omega)$, $\dfrac{1}{p} + \dfrac{1}{q} = 1$ с компактным носителем выполняется равенство

$$\lim_{t \to +0}\|p(tA)f(x) - f(x)\|_{L_p^{-\alpha}(M)} = 0 \qquad (27)$$

для любого компакта $M \Subset \Omega$.

**Доказательство.** Так как пространство $\overset{0}{L}{}_p^s(\Omega)$

$$\|p(tA)u - u\|_{L_p^{-\alpha}(M)} \leq$$

$$\leq \left\|F^{-1}\left[(1+|\xi|^2)^{-\frac{\alpha}{2}} F\left(p(tA)u - u\right)\right]\right\|_{L_p(M)} \leq$$

$$\leq \left\|F^{-1}\left[(1+|\xi|^2)^{-\frac{\alpha}{2}} F\left(p(tA)u - u\right)\right]\right\|_{L_p(M)} \leq$$

$$= \left\|F^{-1}(1+|\xi|^2)^{-\frac{\alpha}{2}} Fp(tA)u - F^{-1}(1+|\xi|^2)^{-\frac{\alpha}{2}} Fu\right\|_{L_p(M)} =$$

$$= \left\|F^{-1}(1+|\xi|^2)^{-\frac{\alpha}{2}-\left[\frac{\alpha-\beta}{2}\right]} FF^{-1}(1+|\xi|^2)^{\left[\frac{\alpha-\beta}{2}\right]} Fp(tA)u - F^{-1}(1+|\xi|^2)^{-\frac{\alpha}{2}-\left[\frac{\alpha-\beta}{2}\right]} FF^{-1}(1+|\xi|^2)^{-\left[\frac{\alpha-\beta}{2}\right]}\right\|$$

. Пусть $s \in R$, $1 \leq p \leq \infty$. Положим

$$\|f\|_p^s = \left\|F^{-1}\left\{(1+|\xi|^2)^{\frac{s}{2}} Ff(\xi)\right\}\right\|_{L_p(R^N)},$$

где $Ff(\xi) = \int_{R^N} \exp(-ix\xi) f(x) dx$ – преобразование Фурье. Пространство Лиувилля $L_p^s(R^N)$ определяются соотношениями

$$L_p^s(R^N) = \left\{ f : f \in S'(R^N), \|f\|_p^s < \infty \right\}.$$

Аналогично, из теоремы 2 следует справедливость следующая

**Теорема 4.** Пусть функция $p(\lambda)$, определенная в $\overline{R_+} = [0, \infty)$, удовлетворяет условиям $p(0) = 1$ и $p(\lambda) \in L_\infty(\overline{R_+}) \cap C([0,\tau])$, $\tau > 0$. Предположим, что $\alpha_0 > \dfrac{N}{p_0}$, $\alpha \geq 0$, $1 < p \leq p_0 \leq 2$ (либо $p = p_0 = 1$), $\varepsilon = N\left(\dfrac{1}{p} - \dfrac{1}{p_0}\right)$, $\beta \geq \alpha_0 + \alpha + \varepsilon$. Тогда для любого распределения $f \in E'(\Omega)$ из класса $\overset{0}{L}{}_q^{-\beta}(\Omega)$, $\dfrac{1}{p} + \dfrac{1}{q} = 1$ с компактным носителем равномерно выполняется равенство $\lim\limits_{t \to +0} p(tA) f(x) = f(x)$ при любом компакта $M \Subset \Omega$.

## § 3.2. О следах средних спектральных разложений, соответствующих самосопряженным эллиптическим псевдодифференциальным операторам, в классах Бесова

Обозначим символом $\overset{0}{B}{}^{\alpha}_{p,q}(\Omega)$ класс функций, принадлежащих пространству Бесова $B^{\alpha}_{p,q}(\Omega)$ и имеющих в $\Omega$ компактный носитель.

Справедлива следующая

**Теорема 3.2.1.** Пусть функция $p(\lambda)$, определенная в $\overline{R}_+ = [0, \infty)$, удовлетворяет условиям

1) $\int\limits_0^\infty |p(\lambda)| \lambda^{\frac{N-\alpha_0}{m}-1} d\lambda < \infty$ и $p(0) = 1$;

2) $p(\lambda) \in C^l(\overline{R}_+)$ и $|p^{(j)}(\lambda)| \leq C_j(1+\lambda)^{-j}$, $j = \overline{0, l}$, где

$l$ – неотрицательное целое число;

3) $l > N\left(\dfrac{1}{2} - \dfrac{1}{p_0}\right)$, $\alpha_0 = \dfrac{N}{p_0}$, $\alpha > 0$, $2 \leq p \leq p_0 < \infty$, $\varepsilon = N\left(\dfrac{1}{p} - \dfrac{1}{p_0}\right)$, $\beta \geq \alpha_0 + \alpha + \varepsilon$ и $1 \leq q < \infty$.

Тогда для любой непрерывной функции $u(x) \in \overset{0}{B}{}^{\beta}_{p\,q}(\Omega)$ с компактным носителем выполняется равенство

$$\lim_{t \to +0} \|p(tA)u(x) - u(x)\|_{B^{\alpha}_{p\,q}(M)} = 0 \qquad (3.2.1)$$

для любого компакта $M \Subset \Omega$.

Доказательство. Пусть $v(x) \in C_0^\infty(\Omega)$. Если символ оператора $A$ имеет вид $\sigma_A = \sigma(y)$, то оператор $p(tA)$ имеет символ $\sigma_{p(tA)} = p(t\sigma(y))$.

Известно, что для любого $k = 0, 1, 2, \ldots$

$$\lim_{t \to +0} \left\| p(tA) D^k v(x) - D^k v(x) \right\|_{L_p(M)} = 0 \ .$$

Используя теоремы вложения, заключаем, что для любого $\varepsilon_1 > 0$ существует $\delta(\varepsilon_1) > 0$, такое, что при $0 < t < \delta$ выполняется неравенство

$$\left\| p(tA) v(x) - v(x) \right\|_{B_{pq}^\alpha(M)} < \varepsilon_1 .$$

Далее, замыканием классов $C^\infty(M)$ в $B_{pq}^\alpha(M)$ является $B_{pq}^\alpha(M)$, т.е. $\overline{C^\infty(M)}\Big|_{B_{pq}^\alpha(M)} = B_{pq}^\alpha(M)$. Поэтому, если $u(x) \in B_{pq}^\alpha(M)$, то для любого $\varepsilon_2 > 0$ существует функция $v(x) \in C^\infty(M)$, для которой

$$\left\| u(x) - v(x) \right\|_{B_{pq}^\alpha(M)} < \varepsilon_2 \ .$$

Для функции $u(x) \in \overset{0}{L}{}_p^{\alpha_0+s+\varepsilon}(\Omega)$ из неравенства (3.1.2) следует справедливость неравенства

$$\left\| p(tA) u(x) \right\|_{L_p^s(M)} \leq \text{const} \cdot \left\| u(x) \right\|_{L_p^{\alpha_0+s+\varepsilon}(R^N)} , \qquad (3.2.2)$$

где $s \geq 0$. Применяя метод вещественной интерполяции, имеем $B_{pq}^\alpha = \left( L_p^{\alpha_1}, L_p^{\alpha_2} \right)_{\theta,q}$, $\alpha_1 \neq \alpha_2$ (см., например [26], [208], [209]).

Следовательно, учитывая (3.2.2), при $s = \alpha'$ и $s = \alpha''$, где $0 \leq \alpha' < \alpha < \alpha''$, получим

$$\left\| p(tA) u(x) \right\|_{B_{pq}^\alpha(M)} \leq \text{const} \cdot \left\| u(x) \right\|_{B_{pq}^{\alpha_0+\alpha+\varepsilon}(R^N)} \qquad (3.2.3)$$

для любой непрерывной функции $u(x) \in \overset{0}{B}{}_{p\,q}^{\alpha_0+\alpha+\varepsilon}(\Omega)$ с компактным носителем, т.е. ограничено равномерно по $t > 0$. Отсюда имеем

$$\|p(tA)(u(x)-v(x))\|_{B_{p\,q}^{\alpha}(M)} \leq C\|u(x)-v(x)\|_{B_{p\,q}^{\alpha_0+\alpha+\varepsilon}(R^N)}.$$

Пусть $v(x) = u_h(x) = \int\limits_{R^N} u(x-hy)\varphi(y)dy$, где $\varphi \in C_0^\infty(R^N)$, $\int\limits_{R^N} \varphi(x)dx = 1$, $\varphi(x) \geq 0$, $\mathrm{supp}\,\varphi = \{x : |x| \leq 1\}$. Если $h$ достаточно мало, то выполнены неравенства

$$\|u(x)-v(x)\|_{B_{p\,q}^{\alpha_0+\alpha+\varepsilon}(R^N)} < \varepsilon_2 \quad \text{и} \quad \|u(x)-v(x)\|_{B_{p\,q}^{\alpha}(M)} < \varepsilon_2.$$

Поэтому в силу неравенства треугольника получим

$$\|p(tA)u(x)-u(x)\|_{B_{p\,q}^{\alpha}(M)} \leq$$

$$\leq C\|u(x)-v(x)\|_{B_{p\,q}^{\alpha_0+\alpha+\varepsilon}(R^N)} + \|p(tA)v(x)-v(x)\|_{B_{p\,q}^{\alpha}(M)} + \|u(x)-v(x)\|_{B_{p\,q}^{\alpha}(M)} <$$

$$< (C+1)\varepsilon_2 + \varepsilon_1 = \varepsilon_3 \quad \text{при} \quad 0 < t < \delta.$$

Выбирая $\varepsilon_2 = \varepsilon_1$, находим, что $\varepsilon_1 = \varepsilon_2 = \dfrac{\varepsilon_3}{C+2}$. Таким образом, для любой непрерывной функции $u(x) \in \overset{0}{B}{}_{p\,q}^{\beta}(\Omega)$ с компактным носителем выполняется равенство

$$\lim_{t \to +0} \|p(tA)u(x)-u(x)\|_{B_{p\,q}^{\alpha}(M)} = 0$$

для любого компакта $M \Subset \Omega$.

Теорема 3.2.1 доказана.

Справедлива следующая

**Теорема 3.2.2.** Пусть функция $p(\lambda)$, определенная в $\overline{R_+} = [0, \infty)$, удовлетворяет условиям $p(0) = 1$ и $p(\lambda) \in L_\infty(\overline{R_+}) \cap C([0, \tau))$, $\tau > 0$. Пусть $\alpha_0 > \dfrac{N}{p_0}$, $\alpha > 0$, $1 < p \le p_0 \le 2$ (либо $p = p_0 = 1$), $\varepsilon = N\left(\dfrac{1}{p} - \dfrac{1}{p_0}\right)$, $1 \le q < \infty$, $\beta \ge \alpha_0 + \alpha + \varepsilon$. Тогда для любой непрерывной функции $u(x) \in \overset{0}{B}{}_{p\,q}^{\beta}(\Omega)$ с компактным носителем выполняется равенство

$$\lim_{t \to +0} \|p(tA)u(x) - u(x)\|_{B_{p\,q}^{\alpha}(M)} = 0$$

для любого компакта $M \Subset \Omega$.

Д о к а з а т е л ь с т в о. Пусть $v(x) \in C_0^\infty(\Omega)$. Если символ оператора $A$ имеет вид $\sigma_A = \sigma(y)$, то оператор $p(tA)$ имеет символ $\sigma_{p(tA)} = p(t\sigma(y))$.

Известно, что для любого $k = 0, 1, 2, \ldots$

$$\lim_{t \to +0} \|p(tA)D^k v(x) - D^k v(x)\|_{L_p(M)} = 0 \ .$$

Используя теоремы вложения, заключаем, что для любого $\varepsilon_1 > 0$ существует $\delta(\varepsilon_1) > 0$, такое, что при $0 < t < \delta$ выполняется неравенство $\|p(tA)v(x) - v(x)\|_{B_{p\,q}^{\alpha}(M)} < \varepsilon_1$.

Далее, замыканием классов $C^\infty(M)$ в $B_{p\,q}^{\alpha}(M)$ является $B_{p\,q}^{\alpha}(M)$, т.е. $\overline{C^\infty(M)}\Big|_{B_{p\,q}^{\alpha}(M)} = B_{p\,q}^{\alpha}(M)$. Поэтому, если $u(x) \in B_{p\,q}^{\alpha}(M)$, то для любого $\varepsilon_2 > 0$ существует функция $v(x) \in C^\infty(M)$, для которой

$$\|u(x) - v(x)\|_{B_{p\,q}^{\alpha}(M)} < \varepsilon_2 \ .$$

Для функции $u(x) \in \overset{0}{L}_p^{\alpha_0+s+\varepsilon}(\Omega)$ из неравенства (3.1.2) следует справедливость неравенства

$$\|p(tA)u(x)\|_{L_p^s(M)} \leq \text{const} \cdot \|u(x)\|_{L_p^{\alpha_0+s+\varepsilon}(R^N)}, \qquad (3.2.4)$$

где $s \geq 0$. Применяя метод вещественной интерполяции, имеем $B_{pq}^{\alpha} = \left(L_p^{\alpha_1}, L_p^{\alpha_2}\right)_{\theta,q}$, $\alpha_1 \neq \alpha_2$ (см., например [26], [208], [209]).

Следовательно, учитывая (3.2.4), при $s = \alpha'$ и $s = \alpha''$, где $0 \leq \alpha' < \alpha < \alpha''$, получим

$$\|p(tA)u(x)\|_{B_{pq}^{\alpha}(M)} \leq \text{const} \cdot \|u(x)\|_{B_{pq}^{\alpha_0+\alpha+\varepsilon}(R^N)} \qquad (3.2.5)$$

для любой непрерывной функции $u(x) \in \overset{0}{B}_{pq}^{\alpha_0+\alpha+\varepsilon}(\Omega)$ с компактным носителем, т.е. ограничено равномерно по $t > 0$. Отсюда имеем

$$\|p(tA)(u(x)-v(x))\|_{B_{pq}^{\alpha}(M)} \leq C\|u(x)-v(x)\|_{B_{pq}^{\alpha_0+\alpha+\varepsilon}(R^N)}.$$

Пусть $v(x) = u_h(x) = \int\limits_{R^N} u(x-hy)\varphi(y)dy$, где $\varphi \in C_0^{\infty}(R^N)$, $\int\limits_{R^N} \varphi(x)dx = 1$, $\varphi(x) \geq 0$, $\text{supp}\,\varphi = \{x : |x| \leq 1\}$. Если $h$ достаточно мало, то выполнены неравенства

$$\|u(x)-v(x)\|_{B_{pq}^{\alpha_0+\alpha+\varepsilon}(R^N)} < \varepsilon_2 \qquad \text{и} \qquad \|u(x)-v(x)\|_{B_{pq}^{\alpha}(M)} < \varepsilon_2.$$

Поэтому в силу неравенства треугольника получим

$$\|p(tA)u(x) - u(x)\|_{B_{pq}^{\alpha}(M)} \leq$$

$$\leq C\|u(x)-v(x)\|_{B_{pq}^{\alpha_0+\alpha+\varepsilon}(R^N)} + \|p(tA)v(x)-v(x)\|_{B_{pq}^{\alpha}(M)} + \|u(x)-v(x)\|_{B_{pq}^{\alpha}(M)} <$$

$$< (C+1)\varepsilon_2 + \varepsilon_1 = \varepsilon_3 \quad \text{при} \quad 0 < t < \delta.$$

Выбирая $\varepsilon_2 = \varepsilon_1$, находим, что $\varepsilon_1 = \varepsilon_2 = \dfrac{\varepsilon_3}{C+2}$. Таким образом, для любой непрерывной функции $u(x) \in \overset{0}{B}{}^{\beta}_{p\,q}(\Omega)$ с компактным носителем выполняется равенство

$$\lim_{t \to +0} \| p(tA)u(x) - u(x) \|_{B^{\alpha}_{p\,q}(M)} = 0$$

для любого компакта $M \Subset \Omega$.

Теорема 3.2.2 доказана.